\documentclass[12pt]{amsart}

\usepackage{tikz}
\usetikzlibrary{matrix}
\usepackage{tikz,tikz-cd}
\usepackage{tikz-cd}

\usepackage{subcaption}  
\usepackage{caption} 

\tikzset{
  symbol/.style={
    draw=none,
    every to/.append style={
      edge node={node [sloped, allow upside down, auto=false]{$#1$}}}
  }
}
\usepackage{parskip}
\usepackage[new]{old-arrows}
\usepackage{amsmath}
\usepackage{hyperref}
\usepackage{amssymb}
\usepackage{latexsym}
\usepackage{amscd}
\usepackage{graphicx} 
\usepackage{amsthm}
\usepackage{mathrsfs}
\usepackage{xypic}
\usepackage{bm}
\usepackage{color}
\usepackage{ulem}

\newdimen\AAdi%
\newbox\AAbo%
\def\AAk#1#2{\s_etbox\AAbo=\hbox{#2}\AAdi=\wd\AAbo\kern#1\AAdi{}}%
\def\AAr#1#2#3{\s_etbox\AAbo=\hbox{#2}\AAdi=\ht\AAbo\raise#1\AAdi\hbox{#3}}%
\font\tenmsb=msbm10 at 12pt \font\sevenmsb=msbm7 at 8pt
\font\fivemsb=msbm5 at 6pt
\newfam\msbfam
\textfont\msbfam=\tenmsb \scriptfont\msbfam=\sevenmsb
\scriptscriptfont\msbfam=\fivemsb
\def\Bbb#1{{\tenmsb\fam\msbfam#1}}
\newtheorem{thm}{Theorem}[section]
\newtheorem{lem}[thm]{Lemma}
\newtheorem{cor}[thm]{Corollary}
\newtheorem{rem}[thm]{Remark}

\newtheorem{defi}[thm]{Definition}

\usepackage{stmaryrd} 
\usepackage{centernot} \usepackage{mathtools}
\usepackage{xcolor}

\newenvironment{pprooff}{\medskip \noindent
{\bf Proof of Theorem \ref{ToO}.}}{\hfill \rule{.5em}{1em}
\\}

\newenvironment{proofff}{\medskip \noindent
{\bf Proof of Theorem \ref{t2}.}}{\hfill \rule{.5em}{1em}
\\}
\newenvironment{prooffff}{\medskip \noindent
{\bf Proof of Theorem \ref{vt2}.}}{\hfill \rule{.5em}{1em}
\\}
\newenvironment{pf}{\medskip \noindent
{\bf Proof of Main Theorem.}}{\hfill \rule{.5em}{1em}
\\}

\newenvironment{pf0}{\medskip \noindent
{\bf Proof}.}{\hfill \rule{.5em}{1em}
\\}
\newenvironment{pfa}{\medskip \noindent
{\bf Proof of Theorem \ref{TA}.}}{\hfill \rule{.5em}{1em}
\\}

\newcommand*{\Lcorner}{%
    \mathchoice%
        {\mathrel{\makebox[7pt][c]{\rule{.4pt}{7.5pt}\rule{5pt}{.4pt}}}}
        {\mathrel{\makebox[7pt][c]{\rule{.4pt}{7.5pt}\rule{5pt}{.4pt}}}}%
        {\mathrel{\makebox[5.5pt][c]{\rule{.4pt}{5.25pt}\rule{3.5pt}{.4pt}}}}%
        {\mathrel{\makebox[4pt][c]{\rule{.4pt}{3.75pt}\rule{2.5pt}{.4pt}}}}%
}

\newcommand{\Section}[2]{\setcounter{equation}{0}
\allowdisplaybreaks
\section[#1]{#2}}

\def\pd#1#2{\frac {\partial #1}{\partial #2}}

\def\a{\alpha}

\def\p#1{\partial #1}

\def\de{\delta}

\def\Om{\Omega}
\def\th{\theta}

\def\R{\Bbb{R}}

\subjclass[2010]{~53A10, ~53A07, ~53C42, ~53C40.}
\begin{document}
\pagenumbering{Roman}\setcounter{page}{1}

\pagenumbering{arabic} \setcounter{page}{1}
\title[On the non-existence of minimal graphs]
 {On the non-existence of solutions of the Dirichlet problem for the minimal surface system}
 \author{Yongsheng Zhang}
\address{Academy for Multidisciplinary Studies, Capital Normal University, Beijing 100048, P. R. China}
\email{yongsheng.chang@gmail.com}

\date{}


\begin{abstract}
In this paper we study the non-existence
of solutions to the Dirichlet problem for minimal graphs of codimension $\geq 2$, 
including certain situations over domain $\Omega$ even with non-$C^1$ boundary $\p \Omega$.
\end{abstract}
\maketitle

\renewcommand{\proofname}{\it Proof.}

\Section{Introduction}{Introduction}\label{S1}
Let $\Omega \subset \mathbb R^{n+1}\ (n\geq1)$ 
be an open bounded 
domain 
with boundary $\p \Omega$ of class $C^r$ for $r\geq 2$.
 {The \textbf{Dirichlet problem for minimal surfaces}} (cf. \cite{j-s,b-d-m, de, m1,l-o})
  asks, for 
 a given boundary (vector-valued) function $f:\partial \Omega \rightarrow \mathbb R^{m+1}$ of class $C^s$ with $0\leq s\leq r$,
         what kind of and how many functions 
         exist in 
         $C^0(\overline{\Omega};\R^{m+1})\bigcap \text{Lip}(\Omega;\R^{m+1})$
        such that each of them,
        $F:x=(x^1,\cdots,x^{n+1})\in \Omega\mapsto F(x)=\big(F^1(x),\cdots,F^{m+1}(x)\big)$, 
        is a weak solution to the minimal surface system
         \begin{equation}\label{ms}
         \left\{\begin{array}{cc}
         \sum\limits_{i=1}^{n+1}\dfrac{\p}{\p x^i}\bigg(\sqrt{g}g^{ij}\bigg)=0, & \ \ j=1,\cdots,n+1,\\
         &\\
         \sum\limits_{i,j=1}^{n+1}\dfrac{\p}{\p x^i}\left(\sqrt{g}g^{ij}\dfrac{\p F^\a}{\p x^j}\right)=0, & \ \ \ \alpha=1,\cdots,m+1,
         \end{array}
         \right.
         \end{equation}
 where $$g_{ij}=\de_{ij}+\sum\limits_{\a=1}^{m+1}\pd{F^\a}{x^i}\pd{F^\a}{x^j},\,\ \ \left(g^{ij}\right)=\big(g_{ij}\big)^{-1}\  \text{ and }\ g=\det\big(g_{ij}\big),$$
 satisfying the Dirichlet boundary condition
         $
          F|_{\partial \Omega}=f.
          $
          
        
        
       When $m=0$ and $\Omega$ has mean-convex $\p\Omega$, 
       according to
        Douglas \cite{d},
                    Rad\'o \cite{r,r2},
                    Jenkins-Serrin \cite{j-s}
                    and Bombieri-de Giorgi-Miranda \cite{b-d-m},
                   with respect to any continuous boundary function $f$ 
                    there exists a unique Lipschitz solution;
                    by de Giorgi \cite{de} and Moser \cite{m1},
                    the solution is in fact analytic in the interior;
                    if furthermore $\p\Omega$ is convex, then the solution has an area-minimizing graph.

When $m\geq 1$, dramatically different behaviors may occur.
        For $\Om=\mathbb D^{n+1}(1)$ (open unit ball of $\mathbb R^{n+1}$)
        with $\mathbb S=\p \mathbb D^{n+1}$,
         the groundbreaking work
          \cite{l-o}
discovered 
     three phenomena  
           about solutions
           for the Dirichlet problem of minimal surface system \eqref{ms}:
           (1) non-uniqueness, (2) irregularity and (3) non-existence,
           with
          suitably chosen boundary condition $\{n, m, f:\mathbb S^n\rightarrow \R^{m+1}\}$ for each of (1), (2) and (3).

          
         To be more precise, 
          Lawson and Osserman established the following  about (3).
                    \begin{thm}[\cite{l-o}]\label{t1}
          Let $\eta:\mathbb S^{n}\rightarrow \mathbb S^m\subset \mathbb R^{m+1}$ be of class $C^2$.
          Suppose $n>m$ and $[\eta]\neq 0\in \ \pi_{n}(S^m)$. 
          Then there exists
           $R_\eta\in\mathbb R_+$ such that 
              when
               constant $R\geq R_\eta$ there is no solution of \eqref{ms}
          for 
          $f=\eta_R
               \triangleq
               R\cdot \eta$.
          \end{thm}
          
         By the implicit functional theorem (see \cite{n}),
         there are always analytic solutions to \eqref{ms} for $f=\eta_R$ whenever $|R|$ is sufficiently small.
         However,
         Theorem \ref{t1}
         indicates
        the  non-existence of solutions to \eqref{ms}
         when $|R|$ is quite large for $\eta$ with $[\eta]\neq0$.
    Therefore, 
    it leads to  the {\bf Lawson-Osserman philosophy} about (2) which suspects that
         there should exist  some $R_0\in (0,R_\eta)$ such that \eqref{ms} has some irregular solution for $f=\eta_{R_0}$.
         Based on the celebrated Hopf maps between spheres,
         Lawson and Osserman 
         constructed
         three examples of cone type irregular  minimal graphs accordingly with suitable slop $R_0$ for each of them.
                  
         In the joint work \cite{x-y-z0} 
          we systematically
          generalized the construction of cone type minimal graphs.
          The boundary function $f=\tan \th \cdot \eta$
       where $\eta$ 
             is an LOMSE 
              composed by a Hopf fibration $\pi$ and a homothetic minimal immersion $\iota$
                             \begin{equation}\label{decLOMSE}\ \ \
        \begin{tikzcd}[/tikz/cells={/tikz/nodes={shape=asymmetrical
  rectangle,text width=1.5cm,text height=2ex,text depth=0.3ex,align=center}}]
 \mathbb S^n  \arrow{d}{\pi} \arrow{dr}{\eta} & \\
  \mathbb KP^l \arrow{r}{\iota} &   \mathbb S^m
  \end{tikzcd}
   \end{equation}
         and $\th$ is the unique acute angle for
              \begin{equation}\label{Sg1}
         \Big\{ \big(\cos\theta\cdot x, \sin\theta\cdot \eta(x)\big):\ x\in \mathbb S^{n}\Big\}
         \end{equation}
          to be minimal in $\mathbb S^{m+n+1}$.
          In such way, 
          we obtained
          (see \S \ref{S4} for details) uncountably many boundary data $\{n, m, f=\eta_{\tan\theta}\}$,
          such that \eqref{ms} has a non-$C^1$  solution
         \begin{equation}\label{Sg00}
         F(y)=\left\{\begin{array}{cc}
|y|\cdot f\Big(\dfrac{y}{|y|}\Big) & y\neq 0,\\
0 & y=0.
\end{array}\right.
\end{equation}
        Some (uncountably many) of these boundary data  can support infinitely many analytic solutions as well,
          which enriches our understanding on  the complexity of (1).
          
       Except the original three in \cite{l-o},
        examples of $\eta$ in \eqref{decLOMSE}
       do not satisfy the assumption of Theorem \ref{t1}.
       In this paper we generalize the non-existence result.

           \begin{thm}\label{t2}
                              If  $\eta: \mathbb S^n\rightarrow \mathbb R^{m+1}$
                              is  $C^2$  
                              with
                              an embedded submanifold image $N^l=\eta(S^n)$
    where $l<n$
                              and $[\eta]\neq 0 \in \pi_n\left(N\right)$,
                              then  there exists $R_\eta\in\mathbb R_+$ such that 
                              when
                              constant $R\geq R_\eta$ there is no solution of \eqref{ms}
          for $f=\eta_R$.
                              \end{thm}  
                           
          
          Every element in the interior part of (graded) moduli spaces of LOMSEs 
          has  an embedded submanifold image,
          so Theorem \ref{t2} applies.
         Actually in Theorem \ref{vt2}
          we can show the same non-existence result
          for
           $\eta=\tilde \iota\circ \pi$ 
           where 
           $\tilde \iota: P \rightarrow \mathbb R^{m+1}$ is an arbitrary immersion (i.e. no requirement on the image),
to include all LOMSEs.
                     As a result
                     the existence of LOMSEs 
                     falls into the  prediction
                     according to the non-existence result here
                     and
                      the Lawson-Osserman philosophy gets further deepened.

           For more general domains besides discs we derive  the following.
          \begin{thm}\label{otherOmega}\label{ToO}
          Let $\Omega\subset \mathbb R^{n+1}$ be a bounded domain with $\p\Omega$ of class $C^2$
          and 
          $\eta:\p\Omega\rightarrow \mathbb R^{m+1}$ be of $C^2$ with  image $N$ of Hausdorff dimension $l<n$.
          If 
          $\eta$
          has
          no continuous extension $\bar \eta$ over $\overline\Omega$ (valued in $\mathbb R^{m+1}$)
      with
           image $\bar \eta(\overline\Omega)$  completely contained in the $\epsilon$-neighborhood of $N$
           for some $\epsilon>0$, 
           then
          there exists $R_\eta\in \mathbb R_+$ such that 
          when constant $R\geq R_\eta$ 
          there is no solution of \eqref{ms} for $f=\eta_R$ on $\p \Omega$.      
          \end{thm}
           
           Examples of $(\Omega, \eta)$ satisfying the requirements can have great diversity,
           e.g.
           $\big(\mathbb D^{n+1}(2)\sim \overline{\mathbb D^{n+1}(1)}, \eta \big)$ with $\eta(S^n(2))\bigcap \eta(S^n(1))=\emptyset$,
           $\big(\mathbb D^n(\frac{1}{2})\times S^1(1),\eta\big)$ with $\eta(x,t)=x$ where $x\in S^{n}(\frac{1}{2})$ and $t\in S^1$,
           and the case of $\Omega\subset \mathbb R^{n+1}$ to be diffeomorphic 
                      to an angular tubular neighborhood $\bold D_\pm$ ($\Subset \mathbb S^{n+1}\sim M_\mp$) 
                        of a focal submanifold $M_\pm$ of an isoparametric foliation of $\mathbb S^{n+1}$
           with $\eta$ induced by the canonical projection of $\p\bold D_\pm\rightarrow M_\mp$
        [One can use the stereographic  projection taking a point of $\mathbb S^{n+1}\sim \big(\bold D_\pm\bigcup M_\mp\big)$
          as a pole
          to realize $\bold D_\pm$, $M_\mp$ and $\eta$ 
          in Euclidean space].
           

           In \S \ref{S53}
          by a generalized version of Lemma \ref{LLO},
          we can get an integral current version of Theorem \ref{ToO}
          and we include some figures to illustrate its applications for domains with non-$C^1$ boundaries.



{\ }

\Section{Preliminaries}{Preliminaries}\label{S2}

In this section we shall recall several lemmas which will be employed in the proofs of our results.
The first is that  if there is a (local) solution around $\p \Omega$ for \eqref{ms} then  it is of class $C^r$ near $\p \Omega$.

\begin{lem}[Allard boundary regularity result in the setting of minimal graphs]\label{TA}
Assume $\p\Omega$ to be of class $C^{1,\alpha}$ for some $\alpha\in(0,1]$
and 
the boundary function $f$ of class $C^{1,\alpha}$.
If $F$ is a Lipschitz solution to the Dirichlet problem \eqref{ms},
then there exists a neighborhood $U$ of $\p \Omega$ such that $F\in C^{1,\alpha}(U\bigcap \overline\Omega)$.
\end{lem}

 \begin{rem}
 A proof will be provided in  Appendix.
 \end{rem}
 
 A  very useful formula about the mass (generalized volume) of a minimal variety 
  (viewed as an integral current of multiplicity one)
 in Euclidean space
  is the following.
 \begin{lem}[e.g. \cite{bl, l-o}]\label{LLO}
 Let $G$ be an oriented $(n+1)$-dimensional minimal variety
  in $\mathbb R^{m+n+2}$ 
 which is a regular $C^1$ submaifold in some neighborhood of $\p G$.
 Then
 \begin{equation}\label{52}
           \bold M(G)=\frac{1}{n+1}\int_{\p G}\big<\nu, p\big>*1_{\p G},
 \end{equation}
           where $p$ is the position vector and $\nu$ the unit exterior normal field to $\p G$.
 \end{lem}

   A well-known property for minimal varieties in Euclidean space is the following.
\begin{lem}[Density Monotonicity, e.g.  \cite{fe,fle,bl, CM, LS}]\label{dmono}
Let $G$ be an oriented $(n+1)$-dimensional minimal variety of $\mathbb R^{m+n+2}$. 
Then the density function 
$$\Theta(G,\tilde p,d):=\dfrac{\bold M(G\bigcap \bold B_d(\tilde p))}{\omega_{n+1}d^{n+1}}$$
where
 $\tilde p\in \mathbb R^{m+n+2}$,
$\bold B_d(\tilde p)$ the solid ball of radius $d$ centered at $\tilde p$
and
$\omega_{n+1}$ the volume of unit ball in $\mathbb R^{n+1}$,
is increasing  for $0< d< \text{dist}(\tilde p, \p G)$.
\end{lem}

           Lawson and Osserman applied
           a contradiction argument 
          to prove  Theorem \ref{t1}.
           Suppose 
           there exist a sequence $\{R_i\}\nearrow +\infty$
such that 
        $F_i$ solves the Dirichlet problem \eqref{ms} for 
$f=\eta_{R_i}$
with graph
denoted  by $G_i$.
           The proof consists of two main parts.
          
           Part I.
           By the topological assumption,
           there exists
           $\tilde p_i=(x_i,0)$ in $G_i$ for each $i$
           and hence
           by Lemma \ref{dmono}
           it follows
       $
 \bold M(G_i)\geq \omega_{n+1}R_i^{n+1}.
    $

            Part II.
            By Lemma \ref{LLO}, it is not hard to see that for $R_i \geq 1$ one has
 \begin{equation}\label{522}
 \aligned
\bold M(G_i)
&\leq \frac{1}{n+1}\int_{p\in \p G_i}\|p\|*1_{\p G_i}(p)
\\
&\leq \frac{\sqrt{1+R_i^2}}{n+1} \, R_i^m \cdot \big(\text{$n$-dimensional volume of the graph of }  \eta \big)
\endaligned
\end{equation}

                                    Thus
       $CR_i^{m+1}  \geq \bold M(G_i)\geq \omega_{n+1}R_i^{n+1}$
    gives a contradiction 
    for
     $R_i\gg 1$  as $n>m$.

   {\ }
           
 \Section{Method for Theorem \ref{t2}}{Method for Theorem \ref{t2}}\label{S3}
           In this section we shall extend Theorem \ref{t1} by understanding the topological condition from a different angle.
         Still  assume that $\Omega=\mathbb D^{n+1}$ at this stage and that
          the image 
          is an embedded submanifold $N^l=\eta(S^n)$ in $\mathbb R^{m+1}$ of dimension $l<n$.    
           
           \begin{defi}
            Let $N$ be an embedded submanifold (for simplicity $\p N=\emptyset$)
            and
            $\epsilon_0(N)$ the largest number for
           the $\epsilon_0(N)$-disk normal bundle of $N$ being diffeomorphic to a neighborhood $\mathscr DN$ of $N$ through the exponential map restricted to normals.
           We call $\epsilon_0(N)$ the normal injectivity radius and
           $\mathscr DN$ the collapsible neighborhood of $N$.            
           \end{defi}
           \begin{rem}
            No essential difference for $N$ with $\p N\neq \emptyset$.
         Given $N\subset \mathbb R^{m+n+2}$ and $R\in \mathbb R_+$, $R\cdot \epsilon_0(N)$ is 
          the normal injectivity radius of $R\cdot N$
   and      $R\cdot\mathscr DN=     \mathscr D(R\cdot N)$.   
           \end{rem}
        Now we are ready to prove Theorem \ref{t2}.
        
\begin{proofff}
 Assume there exist a sequence $R_i \nearrow +\infty$
with $F_i$ being a solution of the Dirichlet problem \eqref{ms} for $f=\eta_{R_i}$
with  graph denoted by $G_i$. 
   
As in the Part II for the proof of Theorem \ref{t1} we can $
           \bold M(G_i)\leq C\cdot R_i^{l+1}.
         $

        For a lower bound of $\bold M(G_i)$, we make the following observation.

          {\bf Claim ($\star$).} For each $i\in \mathbb N$, there exists $\tilde p_i \in G_i$ with $\text{dist}(\tilde p_i, \p G_i)\geq\epsilon_0(N)\cdot R_i$.
          
                   Then Claim ($\star$) and Lemma \ref{dmono} imply that
                                                                 \begin{equation}\label{lobd2}
                                                                 \bold M(G_i)\geq \omega_{n+1}\big(\epsilon_0(N)\big)^{n+1}R_i^{n+1}.
                                                                 \end{equation}
                      Therefore, by  contradiction when $R_i\gg 1$
                      Theorem \ref{t2} gets proved.

               Now let us show Claim ($\star$) by contradiction.
               Assume that, for some $i\in \mathbb N$, $G_i$ is entirely contained in    
                 the $R_i\cdot\epsilon_0(N)$-neighborhood of $\p G_i$. 
               Then the latter $(m+1)$-components of $G_i$ tell that
               $$F_i(\Omega)\subset\mathscr D\big(F_i(\p \Omega)\big)=R_i\cdot \mathscr D(N)$$
               and, in particular,
                          $\big\{F(x,r):=F_i(rx)\ \big|\ x\in \mathbb S^n\text{ and }r\in[0,1]\big\}
                          \subset R_i\cdot \mathscr D(N).
                          $
                 Let $\pi^\perp$ be the retraction of the collapsible neighborhood $\mathscr D(F_i(\p \Omega))$ to $F_i(\p \Omega)$ along normals.
Since 
          $
                \big\{H(x,r):=\frac{1}{R_i} \cdot \big(\pi^\perp\circ F\big)(x,r)\ \big | \ x\in \mathbb S^n\text{ and }r\in[0,1]\big\} \subset N,$      
                 it is clear that $H(\cdot, \cdot)$ establishes a null homotopy for $\eta$ from $\mathbb S^n$ to $N$.
                   This contradicts with $[\eta]\neq0\in \pi_n(N)$ and completes the proof of Claim ($\star$).
\end{proofff}
                    
                    \begin{rem}\label{r22}
                    Decompose $\tilde p_i$ as $(\tilde p'_i,\tilde p''_i)\in \mathbb R^{n+1}\oplus\mathbb R^{m+1}$.
                    Then $p''_i$ can come from the boundary of the collapsible neighborhood. 
                    For $N=\mathbb S^{m}$, considered in \cite{l-o}, 
                    $\tilde p_i$ can always be chosen to be  $(y_i,0)$ for some $y_i\in\Omega$.
Note that  $\tilde p_i''=0\in \p \big(\mathscr D\mathbb S^m\big)\subset \R^{m+1}$.
                    \end{rem}
  
  As what we have seen, the property in the Claim ($\star$) captures a key  ingredient in the proof.
  In fact based on this observation certain situations not having embedded submanifold image
  can also be considered.
  
  {\ }
                    
\Section{Application to LOMSEs}{Application to LOMSEs}\label{S4}       
We shall see that Theorem \ref{t2} applies for most of LOMSEs directly in \S \ref{S41}.
In \S \ref{S42}
we derive Theorem \ref{vt2} to include some situations involving immersed submanifold image.
In particular, the non-existence result can be set up for all LOMSEs.
             
\subsection{Application to interior LOMSEs}\label{S41}

          Recall an LOMSE $\eta$
          in \eqref{decLOMSE} 
         is a composition
          of a Hopf fibration $\pi$ to projective spaces $(\mathbb KP^l,h)$ with
          a minimal isometric immersion $\iota$ of $(\mathbb KP^l,\lambda^2h)$ into Euclidean sphere of unit radius, i.e.,
          $
          \eta=\iota\circ \pi.        
          $
          
%
                                             %
                                             %
                                                Note that there are three families 
                                                $
                                                    \mathbb S^{2s+1}\rightarrow \mathbb CP^s$, 
                                                    $
                                                          \mathbb S^{4s+3}\rightarrow \mathbb HP^s$ 
                                                          and 
                                                                $\mathbb S^{15}\rightarrow \mathbb OP^1$
                                                for the choice of $\pi$,
                                               however 
                                               the (graded) moduli space of $\iota$   in \eqref{decLOMSE} 
 consists of a constellation of polyhedra 
given by suitable frames of the eigenspace of $\tilde k$-th eigenvalue of the Laplacian operator on $(\mathbb KP^l,\lambda^2h)$
(see \cite{c-w, ma1, ma2, oh,u,to,to2}).

                                                The (graded) moduli space of $\iota$ is quite beautiful. The standard isometric immersion $\eta_0$ given by orthonormal basis of eigenfunctions for the $\tilde k$-th eigenvalue lies in
                                              the interior of moduli space.
                                                Since its image $N_0$ is an orbit of corresponding group action (with certain equivariant property), 
                                                $N_0$ is an embedded submanifold in the Euclidean space. 
                                                All other interior point $\eta$ is simply $A\eta_0$,
                                                where $A$ is 
                                                a full rank square matrix $A$ of size of the dimension of the target Euclidean space,
                                                satisfying
                                                \begin{equation}\label{ExA}
                                                A=\sqrt{C+I}
                                                \end{equation}
                                              for symmetric nonnegative matrix $C+I$ with $C\in W_2$.
                                              Here $C\in W_2$ is an equivalent algebraic requirement for 
                                              \begin{equation}\aligned
                                              \|A\eta_0(x)\|&=\|\eta_0(x)\|,\,  \\
                                              \big<A\eta_0(x), Ad\eta_0(v_x)\big>&=\big<\eta_0(x),d\eta_0(v_x)\big>=0,\, \\
                                               \big<Ad\eta_0(v_x), Ad\eta_0(w_x)\big>&=\big<d\eta_0(v_x), d\eta_0(w_x)\big>,
                                               \endaligned
                                              \end{equation}
                                              where $x\in \mathbb S^n$ and $v_x,\, w_x\in T_x \mathbb S^n$ are arbitrary.
                                              The moduli space is described by this compact convex body $W_2$ 
                                              (see \cite{c-w, wa, to, to2} for details).
                                          As a result, each interior point $\eta=A\eta_0$ of the moduli space 
                                          has an embedded submanifold image $N=AN_0$ in an Euclidean sphere.
                                               %

  Moreover, according to \cite{l} and \cite{w-z}, 
$\iota$ is a finite covering on $N$,
  so
  the diagram
        $$
        \begin{tikzcd}[/tikz/cells={/tikz/nodes={shape=asymmetrical
  rectangle,text width=1.5cm,text height=2ex,text depth=0.3ex,align=center}}]
  & \mathbb KP^l  \arrow{d}{\iota}  & \\
 \mathbb  S^n  \arrow{r}{\eta}\arrow{ur}{\pi} &   N^l
  \end{tikzcd}
   $$
 has the homotopy lifting property.
This means $\eta$ is null-homotopic if and only if so is $\pi$.
It is well known that Hopf fibration $\pi$ is not homotopic to a constant map.
          For example, one can use the long exact sequence for homotopy groups,
          called the homotopy sequence of the fibering (see \cite{BT} or \cite{St})
    $$
          \begin{tikzcd}[/tikz/cells={/tikz/nodes={shape=asymmetrical
  rectangle,text width=0.1cm,text height=0ex,text depth=0ex,align=center}}]
         \arrow{r}{\p_*}    &     \pi_n(\mathbb S^{n-l}) \text{\ \ }  &  \arrow{r}{i_*}    & \pi_n(\mathbb S^n)  &  \arrow{r}{\pi_*} & \pi_n(\mathbb KP^l) & \arrow{r}{\p_*} & \pi_{n-1}(\mathbb S^{n-l})
        \cdots
            \end{tikzcd}
          $$
           and 
          apply the celebrated finiteness result due to Serre
         $-$ homotopy groups of spheres are all finite except
          for those of $\pi_s(\mathbb S^s)$ or $\pi_{4s-1}(\mathbb S^{2s})$.
         By the finiteness of $\pi_{2s+1}(\mathbb S^1)$, 
         $\pi_{4s+3}(\mathbb S^3)$ 
         and $\pi_{15}(\mathbb S^7)$, 
          $i_*$ has to be a zero map for
            $$
          \begin{tikzcd}[/tikz/cells={/tikz/nodes={shape=asymmetrical
  rectangle,text width=0.38cm,text height=0ex,text depth=0ex,align=center}}]
         \arrow{r}{\p_*}    &     \pi_{2s+1}(\mathbb S^{1}) \text{\ \ }  &  \arrow{r}{i_*}    & \pi_{2s+1}(\mathbb S^{2s+1})  &  \arrow{r}{\pi_*} & \pi_{2s+1}(\mathbb CP^s) \ \cdots\\
         \arrow{r}{\p_*}    &     \pi_{4s+3}(\mathbb S^{3}) \text{\ \ }  &  \arrow{r}{i_*}    & \pi_{4s+3}(\mathbb S^{4s+3})  &  \arrow{r}{\pi_*} & \pi_{4s+3}(\mathbb HP^s) \ \cdots\\
         \arrow{r}{\p_*}    &   \  \  \pi_{15}(\mathbb S^{7}) \text{\ \ }  &  \arrow{r}{i_*}    & \ \ \pi_{15}(\mathbb S^{15})  &  \arrow{r}{\pi_*} & \ \ \pi_{15}(\mathbb OP^1) \ \, \cdots
            \end{tikzcd}
          $$
respectively.
          Consequently $\pi_*$ is injective. 
          So $[\pi]\neq 0\in \pi_n(\mathbb KP^l)$ and $[\eta]\neq 0 \in  \pi_n(N^l)$.
          Thus, 
          Theorem \ref{t2} applies for all interior LOMSEs in the moduli space.
           
            \begin{cor}\label{c1}
          For any interior LOMSE $\eta: \mathbb S^n\rightarrow \mathbb S^{m}$,
          there exists $R_\eta\in\mathbb R_+$ such that 
          when 
                              constant $R\geq R_\eta$ there is no solution of \eqref{ms}
          for  $f=\eta_R$.
          \end{cor}


          %
        
          
\subsection{Application to boundary LOMSEs}\label{S42}
              In general, the (graded) moduli space of isometric minimal immersions of projective space into Euclidean space 
              may be  subtle in its boundaries, e.g. see \cite{d-z} for some exploration.
              
              It seems still unclear whether the image of a boundary point of the moduli space is always an embedded submanifold or not.
        There might be a possibility that after smashing certain directions (when $A$ in \eqref{ExA} is not of full rank)
              the set $N_0=\eta_0(\mathbb S^n)\subset \mathbb S^m$
              is then immersed into Euclidean spheres $S^{\tilde m}$ (of smaller dimension)
              with self-intersecting points around which local patches do not coincide.
              
              In fact in this subsection, we shall prove a bit more than expected.
              \begin{thm}\label{vt2}
              Let $\eta=\tilde \iota\circ \pi: \mathbb  S^n\rightarrow \mathbb R^{\tilde m+1}$
              where $\pi$ is a Hopf fibration to $\mathbb KP^l$  and $\tilde \iota$ an arbitrary $C^2$ immersion of $\mathbb KP^l$ into $\mathbb R^{\tilde m+1}$.
                           Then  there exists $R_\eta\in\mathbb R_+$ such that for constant $R\geq R_\eta$ there is no solution of \eqref{ms} for $f=\eta_R$.
              \end{thm}
                              \begin{cor}\label{t3}
           For {every} LOMSE $\eta$ in \eqref{decLOMSE},
       there exists $R_\eta\in\mathbb R_+$ such that 
       when 
                              constant $R\geq R_\eta$ there is no solution of \eqref{ms}
          for  $f=\eta_R$.
           \end{cor}
           
            We shall need some concept
             which will be useful in the proof for immersions.
             Note that  
              locally 
              $\tilde \iota$
               is an embedding (for some coordinate chart covering $\{U_j\}_{j=1,\cdots, s}$ of $\mathbb KP^l$) 
              and one can pull back the normal bundle of the immersed  $\mathbb KP^l$ in $\mathbb S^{\tilde m}$.
              \begin{defi}
              Let $(\tilde{\mathscr N}, \mathbb KP, \pi^\perp)$ be the pull-back normal bundle for $\tilde \iota: \mathbb KP^l\rightarrow \mathbb R^{\tilde m+1}$
              and
               $\{\epsilon_i\}$ 
                the  injectivity normal radii
             of $\tilde \iota(U_j)$ in $\mathbb R^{\tilde m+1}$.
               Set $\epsilon_0=\min\{\epsilon_j\}$.
               The $\epsilon_0$-neighborhood $\mathscr D\big(\tilde\iota(\mathbb KP^l)\big)$ 
              is the collapsible neighborhood  for $\tilde \iota$ (subject to $\{U_j\}$).
              \end{defi}
              \begin{rem}
               Let $(\tilde{\mathscr N}_i, \mathbb KP^l, \pi^\perp)$ be the pull-back normal bundle via $\tilde \iota_{R_i}$.
              Then $R_i\epsilon_0$ is the uniform injectivity normal radius for $\tilde \iota_{R_i}$, 
              and $\mathscr D\big(\tilde \iota_{R_i}(\mathbb KP^l)\big)=R_i\cdot\mathscr D\big(\tilde\iota(\mathbb KP^l)\big)$.
              \end{rem}
\begin{prooffff}
          Suppose there exist a sequence $R_i\nearrow +\infty$
supporting solution $F_i$ of \eqref{ms} for each $\eta_{R_i}$
with graph denoted by $G_i$.
Since the rank of $\pi_*$ is $l$, an upper bound on $\bold M(G_i)$ is then given by
 $
           \bold M(G_i)\leq C\cdot R_i^{l+1}.
       $

           For a lower bound, we make the following observation.
           
            {\bf Claim ($\ast$).} For each $i\in \mathbb N$, there exists point $p_i$ in $G_i$ with $\text{dist}(p_i,\p G_i)\geq R_i \epsilon_0$.

          If $(\ast)$ is not true for some $i$, it follows that 
          $G_i\subset (\epsilon_0R_i)\text{-neighborhood of }\p G_i.$
          Consequently, 
          the last $(\tilde m+1)$-component leads to
          \begin{equation}\label{immCase}
          F_i(\mathbb D^{n+1})\subset \mathscr D\big(\tilde \iota_{R_i}(\mathbb KP^l)\big)=
          (R_i \epsilon_0)\text{-neighborhood of }\eta_{R_i}(\mathbb S^n).
          \end{equation}
         We want to show that $F_i$ induces a null-homotopy contradiction of the Hopf fibration.
          
          It is clear that
          $F_i(rx)$ for $x\in \mathbb S^n$ and $r\in [0,1]$
           deforms
          $\eta_{R_i}(\mathbb S^n)$
          to a point in 
          the
          $(R_i \epsilon_0)\text{-neighborhood of }\eta_{R_i}(\mathbb S^n)$.
          Without loss of generality,
          assume every element of the coordinate chart  covering of $\mathbb KP^l$ is homeomorphic to an $l$-dimensional disc.
          With fixed $x\in \mathbb S^n$,
          each curve $F_i(rx)$ for decreasing $r$ in $[0,1]$
          can be projected (through  $\pi^{\perp}$) to 
          involved
          $\tilde \iota_{R_i}(U_j)$
          one patch by another
          tracing the ending fixed curve 
          backward.
          Note that in the overlap
          the projections coincide with each other.
          By the embeddedness of $\tilde \iota_{R_i}|_{U_j}$,
          we get a unique corresponding curve in $\mathbb KP^l$.
          By varying $x\in \mathbb S^n$
          these curves assemble a continuous deformation
            \footnote{The homotopy may not be entirely $\mathbb S^1$-invariant.} 
           (parametrized by $r$ within the target space $\mathbb KP^l$) of the Hopf fibration $\pi$ to a constant map.
          But  $\pi$ is not null-homotopic.
          Contradiction!
              Hence Claim $(\ast)$ is true.

                    Therefore $CR_i^{l+1}  \geq \bold M(G_i)\geq \omega_{n+1}\epsilon_0^{n+1}R_i^{n+1}$
                    which is impossible for $R_i\gg 1$.
 \end{prooffff}

      

   
   {\ }
   
  \Section{Other regions}{Other regions}\label{S5}
  All previous discussions are devoted to the situation with $\Omega=\mathbb D^{n+1}(1)\subset \mathbb R^{n+1}$.
Based on \S \ref{S2}, we can derive the non-existence result for more general  domains.
Let us first give a quick proof of Theorem \ref{ToO}
and then move on its applications to several interesting domain $\Omega$.

\begin{pprooff}
The assumption about $\epsilon$ in Theorem \ref{ToO} is directly the key property that we derive for all previous results
and it indicates a lower bound $C\cdot R^l$ on mass for the graph of solution to \eqref{ms} for $f=\eta_R$.
The upper bound on mass can be derived in the same. 
As $l<n$ we  can prove the non-existence in the statement as before.
\end{pprooff}

All the remaining results in the paper are somehow corollaries of Theorem \ref{ToO}.
To see different reasons to fulfill the assumption in Theorem \ref{ToO}, we elaborate more.
          
   \subsection{Application to region with disconnected $\p\Omega$}\label{S51}
  Now we study the Dirichlet problem \eqref{ms} over region $\Omega=\mathbb D(2)\sim \overline{\mathbb D(1)}$ 
   with $C^2$ boundary functions $\eta_1: S^{n}(1)\rightarrow \mathbb R^{m+1}$ and $\eta_2: S^{n}(2)\rightarrow \mathbb R^{m+1}$.

   \begin{cor}\label{Sepa}
   Assume that $\eta_1(S^n(1))$ and $\eta_2(S^n(2))$ are disjoint sets of Hausdorff dimensions strictly less than $n$.
   Then there exists $R_{\eta_1,\eta_2}\in \mathbb R_+$ such that for constant $R\geq R_{\eta_1,\eta_2}$
   the Dirichlet problem has no solution  for boundary data $\{R\cdot \eta_1$, $R\cdot \eta_2\}$.
   \end{cor}
   \begin{pf0}
         Let $d=\text{dist}\big( \eta_1(S^n(1)),\ \eta_2(S^n(2))\big)>0$.
     Then $\epsilon=\frac{d}{2}$ fulfills the assumption of Theorem \ref{ToO}.
   \end{pf0}
   \begin{rem}\label{notall}
   The result holds for the case $m=0$.
   When $\{\eta_1\equiv 1, \eta_2\equiv 0\}$,
   it says that there are no ``tall catenoids" if two ends are relatively far away from each other.
     \end{rem}
{\setlength{\parindent}{0cm}
   {\bf Example 1.} Let $\eta_1:S^3(1)\rightarrow S^2(1)$ be the first Hopf map and $\eta_2:S^3(2)\rightarrow 2\cdot S^2(1)$ by $x\mapsto 2\cdot \eta_1(\frac{x}{\|x\|})$.
   Then for large $R$, the Dirichlet problem \eqref{ms} for $\{R\cdot \eta_1,\ R\cdot \eta_2\}$ has no solutions.
}

   How about $\eta_1$ being the first Hopf map and $\eta_2(x)=T\circ \eta_1\big(\frac{x}{\|x\|}\big)$? Here 
   \begin{equation*}
T(z, x_3)=\left\{\begin{array}{cc}
\Big(\dfrac{z^2}{\|z\|}, x_3\Big) &\quad \text{\ \ \ where } z\in \mathbb C\sim \{0\} \text{ and } (z,x_3)\in S^2(1)\subset \mathbb R^3,\\
&\\
\big(0, x_3\big) &  \text{ when } z=0 \text{ and } x_3=\pm 1.
\end{array}\right.
\end{equation*}
   Note that, as a self-map of $S^2$, $T$ is smooth.
   Although
   $ \eta_1(S^n(1))= \eta_2(S^n(2))$,
   by the non-null-homotopy of $T$
the same non-existence holds for such situation.
   
   \begin{cor}\label{D12}
   Suppose $n>m+1$. 
  If $C^2$ mappings $\eta_i: S^n(i)\rightarrow \mathbb R^{m+1}\sim \{0\}$ for $i=1,2$ are not homotopic to each other
  and have images of strictly smaller dimensions,
   then there exists $R_{\eta_1,\eta_2}\in \mathbb R_+$ such that when constant $R\geq R_{\eta_1,\eta_2}$
  the Dirichlet problem \eqref{ms} has no solution for $\{R\cdot \eta_1,\ R\cdot \eta_2\}$.
   \end{cor}
   
   \begin{pf0}
   The punch line is that if there existed a solution $F$ for $\{R\cdot \eta_1, R\cdot \eta_2\}$,
   there would be some point $x_R\in \mathbb D(2)\sim \overline{\mathbb D(1)}$ with $F(x_R)=0$
   which can further lead to a contradiction on the mass of the graph of $F$ for large $R$ as in Lawson-Osserman's proof of Theorem \ref{t1} in \S \ref{S2}.
   \end{pf0}
   \begin{rem}
   There is also a corresponding version of the non-existence for the boundary situation $\{\eta_1, R\cdot \eta_2\}$ when $R$ is quite large.
   (cf. Theorem \ref{Sepa}).
   \end{rem}
   Let us briefly mention a more complicated situation inspired by the proof of Theorem \ref{t1}. 
   Denote $\mathbb D_i=\mathbb D(3-i)$ for $i=1,2$.
   Let $\tilde\Omega=\left(\mathbb D_1\sim \overline{\mathbb D_2}\right)\sim\bigcup_{j=3}^k \overline{\mathbb D^{n+1}_j}$
   where $\left\{{\mathbb D_j}\right\}_{j=3,\cdots, k}\Subset \mathbb D(2)\sim \overline{\mathbb D(1)}$ are disjoint
   open balls 
   and $\eta=\{\eta_1,\eta_2,\cdots, \eta_k\}$ the set of $C^2$ boundary functions from $\p\mathbb D_j$ to $\mathbb R^{m+1}\sim \{0\}$ for $j=1,2,\cdots, k$.
   When $n>1$ we can do the connect sum among $\p\mathbb D_j$ for $2\leq j\leq k$ 
   using thin necks
   around disjoint curves (to be decided) in $\tilde \Om$.
   Apart from that, we can always abstractly define a trivial connect sum among  maps $\{\eta_{j}\}_{2\leq j\leq k}$.
   For simple illustration, without loss of generality
   suppose $\eta_j(\p \mathbb D_j)\subset \mathbb S^m$
   and $[\eta_2], [\eta_3]$ are nonzero in $\pi_n(\mathbb S^m)$.
   Take $q_j\in \p \mathbb D_j$ for $j=2,3$
   with $\eta_2(q_2)=\eta_3(q_3)=p\in \mathbb S^m$
   and one can deform $\eta_2,\eta_3$ a bit 
         to make 
         values of $\eta_2$ and $\eta_2$ being constant $p$
         in a small neighborhoods of $q_2$ and $q_3$ respectively.
         Now we can have a trivial connect sum  $\eta_2\#\eta_3$ by defining constant $p$ on the connecting tube region
         and its homotopy class is denoted by $[\eta_2]+[\eta_3]$.

   \begin{cor}
   Suppose that every $\eta_j(\p \mathbb D_j)$ has Hausdorff dimension less than $n$ (which is automatic if $m+1<n$).
   Assume that $m>1$ and that 
   \begin{equation}\label{Hoetaj}
\sum_{j=2}^k[\eta_j]\neq [\eta_1] \in \pi_n(\mathbb R^{m+1}\sim \{0\})\cong\pi_n(\mathbb S^m).
   \end{equation}
   Then there exists $R_{\eta_1,\cdots,\eta_k}\in \mathbb R_+$ such that for constant $R\geq R_{\eta_1,\cdots,\eta_k}$
  the Dirichlet problem \eqref{ms} has no solution  for $\{R\cdot \eta_1,\cdots, R\cdot \eta_k\}$.
\end{cor}
\begin{pf0}
Note that for $\{R\cdot \eta_1,\cdots, R\cdot \eta_k\}$, if there exists a solution  $F$ of \eqref{ms},
then there must be some point $x_R\in \tilde \Omega$ such that $F(x_R)=0$.
Otherwise, one can concretely consider set $\p\mathbb D_2\#_{\gamma_2}\p\mathbb D_3\#_{\gamma_3}\cdots \#\p_{\gamma_{k-1}}\mathbb D_k\subset \tilde \Omega$ 
by thin necks
 around mutually disjoint connecting curves $\gamma_2,\cdots, \gamma_{k-1}$.
As $F$ is fixed for each $R$, here we can take the necks sufficiently thin such that in the sense of homotopy group
the restriction of $F$ on the connect sum set decides the same class
$\sum_{j=2}^k[\eta_j]$ 
by the trivial connect sum among $\{\eta_j\}_{j=2}^k$ which we explain abstractly in the above. 
Therefore, it contradicts with \eqref{Hoetaj}.
Hence we can apply Lemma \ref{dmono} centered at $(x_R,0)$ for a lower bound of the graph of $F$ 
as $R\epsilon=R\cdot\text{dist}\big(0, \bigcup_{j=1}^k \eta_j(\p \mathbb D_j)\big)$.
Meanwhile, the dimension assumption assures the desired  upper bound of mass of the graph of $F$.
Thus a combination of both implies the conclusion.
\end{pf0}
\begin{rem}
Although $\p\mathbb D_2\#_{\gamma_2}\p\mathbb D_3\#_{\gamma_3}\cdots \#\p_{\gamma_{k-1}}\mathbb D_k\cong \p\mathbb D_1$, 
the maps on the necks are a priori unknown.
So one cannot directly apply Theorem \ref{D12}.
\end{rem}
  
  {\ }
  
\subsection{Some other domains and obstruction arising from homologies}\label{S52}
In this subsection we would provide more examples.
Let 
           $$S^1(1)=\left\{(0,\cdots, 0, x_n, x_{n+1}):x_n^2+x_{n+1}^2=1\right\}\subset \mathbb R^{n+1}.$$
Then at each $t\in S^1(1)$ there exists a unique perpendicular normal hyperplane $t^\perp\cong \mathbb R^{n}$
and thus one can write small distance tubular neighborhood of $S^1(1)$ as a product.

                \begin{cor}\label{T521}
                Let $\bold T^{n+1}=\mathbb D^n\big(\frac{1}{2}\big)\times S^1(1)\subset \mathbb R^{n+1}$.
                If a $C^2$ mapping $\eta:\p \bold T\rightarrow \mathbb R^{m+1}$
                has a simply connected embedded submanifold image $N$ of lower dimension  
                and $\eta$ is not homotopic to a constant map between $\p \bold T$ and $N$,
                then there exists $R_\eta\in\mathbb R_+$ such that
                for constant $R\geq R_\eta$ the Dirichlet problem \eqref{ms} on $\bold T$ with boundary data $f=R\cdot \eta$
                has no solutions.
                \end{cor}
                \begin{pf0}
                Note that a solution $F$ for $\eta$
                gives a homotopy between $\eta:\p \bold T\rightarrow N$ and 
                $\eta_0:\p \bold T\rightarrow \pi^\perp \circ F(\{0\}\times S^1(1))\subset N$
                if $F(\bold T)$ is contained in the collapsible neighborhood $\mathscr D(N)$ of $N$.
                However this leads to a contradiction of the simple connectedness of $N$ with the topological assumption on $\eta$.
                So, by choosing $\epsilon$ to be the normal injectivity radius $\epsilon_0(N)$ of $N$, Theorem \ref{ToO} implies the corollary.
                \end{pf0}

{\setlength{\parindent}{0cm}
          {\bf Example 2.}
          Let $m=n-1\geq 2$ and $\eta:\p \bold T\rightarrow S^{n-1}\big(\frac{1}{2}\big)$ by $(x,t)\mapsto x$. 
          Then 
      \begin{equation}\label{map}
                  S^{n-1}\left(\frac{1}{2}\right)
                  \cong 
                  S^{n-1}\left(\frac{1}{2}\right)\times \big\{t_0\big\}
          \longhookrightarrow \p \bold T=S^{n-1}\left(\frac{1}{2}\right)\times S^1\big(1\big)\joinrel\xrightarrow{ \ \ \eta\  \ } S^{n-1}\left(\frac{1}{2}\right)
  \end{equation}
          indicates that $\eta$ is not null-homotopic.
          By Corollary \ref{T521} there is no solution for $f=R\cdot \eta$ over $\bold T$ for large $R$.}

\begin{rem}\label{rkn-1=1}
In Example 2, 
           the result also works when $m=n-1= 1$.
           Note that
           now $N$ is a circle not simply connected.
           However, by the non-null-homotopy of $\eta$ in \eqref{map}, for every solution $F$ to \eqref{ms} for $\eta$ and $t_0\in S^1(1)$,
           the argument in the proof of Theorem \ref{t1}
          indicates the existence of
           $x_{F,t_0}\in \mathbb D^{n+1}(\frac{1}{2})\times \{t_0\}$
           with $F(x_{F,t_0})=0$.
           Hence, Theorem \ref{ToO} with $\epsilon=\frac{1}{2}$ applies.
           \end{rem}
  By the observation in Remark \ref{rkn-1=1}, a variation of Corollary \ref{T521} forms the following.

{\setlength{\parindent}{0cm}
{\bf Example 3.} Given any simple closed curve $\gamma$ in $\mathbb R^3$.
Then its small distance tubular neighborhood $\subset \mathscr D(\text{Image of }\gamma)$ is an embedding submanifold which is diffeomorphic to $\mathbb D^2\times S^1$.
With $\eta$ as in Example 2,
the same type non-existence result holds.
}

          Assuming everything is at least $C^2$, here come more corollaries of Theorem \ref{ToO}.
         
         \begin{cor}
         For $n-k\geq k+1\geq 2$ and domain $\Omega=\mathbb D^{n-k+1}\left(\frac{1}{2}\right)\times S^k(1)$, 
          if  boundary mapping $\eta$ 
          has an embedded submanifold image $N^l$ in $\mathbb R^{m+1}$ with $l<n$, 
          $\eta$ is not null-homotopic between $\p \Omega$ and $N$,
           and $\pi_k(N)=0$, 
          then the Dirichlet problem \eqref{ms} has no solutions for boundary data $f=\eta_R$ when constant $|R|$ is large.
        \end{cor}
        
        In general the obstruction can come from homologies instead.
        
 \begin{cor}\label{T522}          
         Suppose that $k\in \mathbb Z_+$,
          that $(M_1^{n-k+1},\p M_1)$ is
          an embedded manifold with boundary in $\mathbb D^{n-k+1}(\delta)$,
          which can retract to its submanifold $K_1$ of dimension $d<n-2k$,
          and that 
          $M_2^k$ is a closed oriented hypersurface of $\mathbb R^{k+1}$ with normal injectivity radius  greater than $\delta$.
          With the splitting $\mathbb R^{n+1}\cong \mathbb R^{n-k}\times \mathbb R^{k+1}$
          (pointwise, cf.  the beginning of this subsection for $S^2(1)$),
          consider $\Omega=M_1\times M_2\subset \mathbb R^{n+1}$
          and take $\eta:\p M_1\times M_2\rightarrow \p M_1\subset \mathbb R^{n-k+1}$ 
          so that 
          \begin{equation}\label{hol}
          \eta_*\neq 0: H_{n-k}(\p M_1\times M_2;\, \mathcal R       )\longrightarrow H_{n-k}({\p M_1};\, \mathcal R    )\ \text{ where $\mathcal R=\mathbb Z$ or $\mathbb Z_2$}.
          \end{equation}
          Then there exists $R_\eta\in\mathbb R_+$ such that
               when constant $R\geq R_\eta$ the Dirichlet problem \eqref{ms} 
                for boundary data $f=R\cdot \eta$
                has no solutions.
 \end{cor}   
                             \begin{rem}
                             For example, $M_1$ can be a collapsible  neighborhood of an embedded $\mathbb R P^2$ in $\mathbb D^{2026}(\delta)$.
                           As
                        the normal bundle     along $\{0\}\times M_2$
                        in  $\mathbb R^{n-k}\times \mathbb R^{k+1}$
                        is trivial
                        of rank $n-k+1$,
                        we can define an $(n+1)$-dimensional region $\Om$ 
                        by $M_1\times M_2$.
                        Strictly speaking,
                        when the global trivialization of the normal bundle
                        is different,
                        geometrically the region $\Om$ changes accordingly.
                        However, we can use any trivialization as essentially the topology plays a decisive role here.
                             \end{rem}
 \begin{pf0}      
          Given a solution $F$ of the Dirichlet problem for $\eta$,
          if $F(M_1\times M_2)$ is entirely contained in the collapsible neighborhood $\mathscr D(\p M_1)$ of $\p M_1$ in $\mathbb R^{n-k+1}$,
          then 
          $\eta\simeq \eta'$ where $\eta': \p M_1\times M_2\rightarrow \pi^\perp\circ F\left(K_1^d\times M_2^k\right)\subset \p M_1$
          (of dimension $n-k$)
          contradicts with $\eta_*\neq 0$ \eqref{hol} due to the assumption $d<n-2k$.
        Hence, by choosing $\epsilon$ to be the normal injectivity radius $\epsilon_0(\p M_1)$ of $\p M_1$, Theorem \ref{ToO} implies the corollary.          \end{pf0}

   \subsection{Situations associated to isoparametric foliations of spheres}\label{S53}
  Among others, some  explicit examples with obstruction from homologies involving 
   more complicated topologies for the non-existence phenomena
   come from isoparametric foliation of Euclidean spheres.
           Let us recall some very basic knowledge of the isoparametric theory of spheres.

    A closed embedded hypersurface $M$ in 
a unit sphere $\mathbb S^L$ is called isoparametric, by E. Cartan, if it has constant principal curvatures.
             When $M$ is isoparametric, so are its parallels.
             In such way, a foliation of hypersurfaces appears with two exceptional leaves $M_\pm$ $-$ focal submanifolds of higher codimensions.
       
       Assume $\xi$ to be a unit normal vector field along $M$ towards $M_+$,
               $g$ the number of distinct principal curvatures of $M$,
          $\cot \theta_\alpha$ $(\alpha=1,\cdots, g; 0<\theta_1<\cdots<\theta_g<\pi)$
           the principal curvatures with respect to $\xi$
           and
           $m_\alpha$ the multiplicity of $\cot\theta_\alpha$.
           An elegant result \cite{M} by M\"{u}nzner  says that $g$ must be $1,2,3,4$ or $6$
           and $m_\alpha = m_{\alpha+2}$ (indices mod $g$).
           $M_\pm$ have codimensions $m_1+1$ and $m_2+1$ respectively,
           in the way that
   each hypersurface leaf $M$ is 
   the boundary of a $\mathbb D^{m_1+1}$-bundle $\bold D_+$ over $M_+$ 
   and simultaneously also that of a $\mathbb D^{m_2+1}$-bundle $\bold D_-$ over $M_-$.
   Namely, $M=\p \bold D_+=\p \bold D_-$ can be written as a sphere-bundle over $M_\pm$, with projection maps $\pi_\pm$ respectively.

     Let us focus on the case with $g=4$. 
     By the beautiful results of Cecil-Chi-Jensen \cite{CCJ} and Chi \cite{Chi}, 
     a foliation with $g=4$ must be either of OT-FKM type or homogeneous with $(m_1,m_2)=(2,2)$ or $(4,5)$.
     In particular, there are infinitely many examples of isoparametric foliations with unequal $m_1, m_2$ and interesting $M_+$
     and $M_-$ (some could be non-orientable).
     Fix a pair of embeddings $\iota_\pm$ of $\bold D_\pm \hookrightarrow  \mathbb R^L$.
     Of course, there is a huge flexibility to do so, for example one can take suitable stereographic projections.
    As a result, through the diffeomorphism $(\iota_+)^{-1}$ 
    we can identify $\iota_+(M_+)$ with $M_+$
    and the same for $\iota_-(M_-)$ with $M_-$ under $(\iota_-)^{-1}$. 
     As we just need diffeomorphisms for $\bold D_-$ and $M_+$ 
                                  (or $\bold D_+$ and $M_-$) 
                                  in Euclidean spaces respectively,
                                  one can take one stereographic projection 
                                  to embed both $\bold D_-$ and $M_+$ in one Euclidean space $\mathbb R^L$
                                  as mentioned in \S \ref{S1}
                                  and the total space for the minimal graph Dirichlet problem \eqref{ms}
                                  now is $\mathbb R^L\oplus \mathbb R^L$.
                                  Of course, in general 
                                  one can replace the second copy of $\mathbb R^L$ by Euclidean space of less dimension
                                  as what we need is just to embed $M_+$ into an Euclidean space.
     
     \begin{cor}\label{isoparametric}
     Suppose $\text{dim}(M_+)>\text{dim}(M_-)$ of an isoparametric foliation of $\mathbb S^L$.
    Take $\Omega=\iota_-(\bold D_-)\subset \mathbb R^L$ and 
    $\eta=\iota_+\circ\pi_+\circ (\iota_-)^{-1}:\p \Omega=\iota_-(\p \bold D_-)\rightarrow \iota_+(M_+)\subset 
    \mathbb R^L$.
     Then there exists $R_\eta\in\mathbb R_+$ such that
                for constant $R\geq R_\eta$ the Dirichlet problem \eqref{ms} on $\Omega$ with boundary data $R\cdot \eta$
                has no solutions.
     \end{cor}
     \begin{rem}
     By the construction of symmetric Clifford systems for OT-FKM type isoparametric foliations with $g=4$ on spheres,
     most of the time 
     we have $m_2>m_1$
     and thus $\text{dim}(M_+)>\text{dim}(M_-)$.
     \end{rem}
    \begin{pf0}
    Recall that M\"{u}nzner \cite{M} applied the Mayer-Vietoris sequence for
    \begin{equation}\label{plusholo}
    H^q(M;\, {\mathcal R})=H^q(M_+;\, {\mathcal R})\oplus H^q(M_-;\, {\mathcal R})
    \end{equation}
    with injective homomorphisms
    \begin{equation}\label{Holo0}
    \pi_\pm^*: H^{q}(M_\pm;\, {\mathcal R}) \longrightarrow H^q(M;\, {\mathcal R}).
    \end{equation}
    In particular, as $\p \bold D_-=M$, 
    for we get
     \begin{equation}\label{Holo}
     \big(\iota_-\big)^*\circ\eta^*\circ \big((\iota_+)^{-1}\big)^*=\pi_+^*\neq 0: H^{m_1+2m_2}(M_+;\, {\mathcal R})\cong{\mathcal R} \longrightarrow H^{m_1+2m_2}(M;\, {\mathcal R}).
     \end{equation}
    Here we use ${\mathcal R}=\mathbb Z$ if both $M_\pm$ are orientable, and  ${\mathcal R}=\mathbb Z_2$ otherwise.
    
    Let $F$ be any continuous extension of $\eta$ from $\p \Omega $ to $\overline{\Omega}$ with values in $\mathbb R^L$.
          If $F(\Omega)$ is entirely contained in the collapsible neighborhood of $\iota_+(M_+)$,
          then 
          $\eta\simeq \eta'$ where
          $\eta':  \p \Om=\iota_-(\p \bold D_-) \rightarrow \pi^\perp\circ F\big(\iota_-(M_-)\big)\subset \iota_+(M_+)$.
          Here 
          \begin{equation}\label{eta'}
          \eta'= \pi^\perp\circ F\circ\pi_-\circ (\iota_-)^{-1}
           \end{equation}
          and the homotopy between $\eta$ and $\eta'$ can be given as follows.
          Let $ \bold D_-$ be a $\delta$-neighborhood of $M_-$ in $\mathbb S^L$
          which
          is
           diffeomorphic to the $\delta$-ball normal bundle $\mathscr B(\delta)$ of $M_-$ in  $\mathbb S^L$ via the exponential map $\exp^\perp$ restricted to normals.
       Hence 
       we can identify 
$\iota_-( \bold D_-)$ with   $\mathscr B(\delta)=\{(x,v)\,:\, x\in M_- \text{ and } v\in N_x(M_-)\text{ with $\|v\|=\delta$}\}$ 
(and thus omit $(\iota_-)^{-1}$).
  The homotopy 
              $H: \p \Om \times [0,1] \rightarrow \iota_+(M_+)$
                    is given by 
      \begin{equation}\label{eta'homotopy}
               H\big((x,v), t\big)
               =(\pi^\perp\circ F)\big(\iota_-(\exp^\perp_x(tv))\big)
                        \end{equation}
    which is $\eta$ for $t=1$
    and produces $\eta'$ when $t=0$.   
    Therefore, $\eta^*=\eta'^*$.
    But $\pi^\perp\circ F$ in $\eta'$ is from 
    $\iota_-(M_-)$ to $\iota_+(M_+)$.
       By  the dimension assumption $\eta^*=\eta'^*=0$ between cohomology groups of degree $m_1+2m_2$
          contradicting with \eqref{Holo}.
          
Hence, by choosing $\epsilon$ to be the normal injectivity radius $\epsilon_0\big(\iota_+(M_+)\big)$, Theorem \ref{ToO} implies the corollary.   
    \end{pf0}
    \begin{rem}
    In fact, the assumption that $\text{dim}(M_+)>\text{dim}(M_-)$ can  be removed 
               according to 
                        the decomposition \eqref{plusholo} and 
                        the injectivity \eqref{Holo0}.
    Namely, one can always get that
                        $F(\Omega)$ can never be entirely included in the collapsible neighborhood of $\iota_+(M_+)$
                        regardless the relation between 
                        $\text{dim}(M_+)$ and $\text{dim}(M_-)$.
                        Otherwise, 
                        instead consider
                        the exact sequence
                        $$
                   0   \longrightarrow    
                   H^q(
                   M_+;\, {\mathcal R})
                   \xlongrightarrow    
                   {
                   \pi_+^*}
                    H^q(M;\, {\mathcal R})=H^q(M_+;\, {\mathcal R})\oplus H^q(M_-;\, {\mathcal R})
                    \xlongrightarrow    
                    {
                    i^*
                    }
                    H^q(
                    M_-;\, {\mathcal R})
                    \longrightarrow    0
                        $$
                        where 
                        $i$ is the inclusion of $M_-$ in $\bold D_-$.
                        As $\eta'$ factors through $\pi_-$,
                        so the image of $\iota_-^*\circ \eta^*$ equals that of $\iota_-^*\circ\eta'^*$
                        which belongs to both 
                        $
                        H^q(M_+;\, {\mathcal R})
                        $
                        and 
                        $
                      H^q(M_-;\, {\mathcal R})
                        $
                        and hence must be the zero class of $ H^q(M;\, {\mathcal R})$.
                        However this contradicts with \eqref{Holo}.

    Moreover, analogously  to Theorem \ref{vt2},
    we can take $\eta$ to be $\tilde \iota\circ\pi_+\circ(\iota_-)^{-1}$
     where $\tilde \iota$ is an arbitrary immersion of $M_+$ into any Euclidean space $\mathbb R^{m+1}$
    for the same type result.
        \end{rem}

   \subsection{Results with non-$C^1$ domains}\label{S53}
   
   In this subsection we will consider the situation of domains with singularities in boundaries.
   A key point of being able to extend previous results to  this setting is a generalized version of Lemma \ref{LLO}
   valid for every stationary integral current with boundary in Euclidean space.
   
                     Before stating the result,
                     let us prove a useful lemma.

                        \begin{lem}\label{Gl}
                       Let $G$ be a stationary $n$-current in $\R^N$ with boundary $\p G$ , namely it is stationary for variations away from $\p G$.
                       Then
                        \begin{equation}\label{52N0}
           \bold M(G)=\frac{1}{n+1}\int_{Z}\big<\nu, p\big>d\sigma,
 \end{equation}
                         where $p$ is the position vector,
                 $Z$ the generalized boundary,
                 $\sigma$ the generalized boundary measure of $G$
                  and $\nu$ the generalized unit co-normal of $G$.

   \end{lem}
                                      \begin{proof}
                                        Every integral current  in Euclidean space can  be automatically identified
                     as
             an
                     integer multiplicity stationary varifold.
                         We can employ the first variation formula, 2.10 of Chapter 8 in \cite{LS},
                     \begin{equation}\label{1stV0}
                                       \delta G (X)
                                          =
                                          -\int_{\Om\times \R^{m+1}}       \underline H       \cdot      X      \, d\mu_G
                                          -\int_Z \nu\cdot X  \, d\sigma
                    \end{equation}
                           where $Z, \sigma, \nu$
                           are mentioned in the the statement of the lemma,
                           and moreover, 
                           $d\mu_G$ the supported measure by $G$ (of integer multiplicity with respect to the Hausdorff $n$-measure on the support set of $G$),
                           $\underline H$ the generalized mean curvature.
                           Since $G$ is stationary away from the boundary,
                           it follows that the first term on the right hand side of \eqref{1stV0}
                           vanishes.
                           Therefore,
                            \begin{equation}\label{1stV1}
                           \delta G (X)
                                          =
                                          -\int_Z \nu\cdot X  \, d\sigma.
                             \end{equation}
                           Consider the homothety family $\{tG\}$ for $t\in (0,1]$
                           for which on each ``slice" $tG$ the deformed vector field is exactly the restriction $X_t$ of Euler vector of $\R^N$. 
                           Note that  one can have \eqref{1stV1} for each $tG$.
                           Hence, we get
                              \begin{equation}
                          \frac{d\,  \bold M (G_t)}{dt} 
                                          =
                                          -\int_{tZ} \nu_t\cdot X_t  \, d\sigma_t
                                           =-t^n \int_{p\in Z}\nu\cdot p\, d\sigma.
                   \end{equation}
                 Therefore,
                 $$
                                             \bold M (G_t)         =-\frac{1}{n+1}
                                                 \int_{p\in Z}\nu\cdot p\, d\sigma
                 $$
                 and we finish the proof.
\end{proof}

                    Suppose that $\Om$ is an $(n+1)$-dimensional domain 
                    and 
                    that
                                 the graph of  boundary $f$ over $\p \Om$ in $\R^{n+m+2}$
                                 is an integral current of multiplicity one (thus $\p \Om$ needs to be an integral current).
                                 Then every solution $F$ to the minimal surface system \eqref{ms}
                                 has a graph $G$
                                 which
                                 is both inner and outer critical.
                                 Hence $G$ is an  integral current of multiplicity one which is stationary under variations away from its boundary.

   \begin{cor}\label{Gl2}
                       Let $F$ be a (Lipschitz) solution of minimal surface system \eqref{ms}.
                       Suppose that 
                        the graph of  boundary $f$ over $\p \Om$
                        is an integral $n$-current.
                       Then the graph $G$ of $F$ is an $(n+1)$-integral current,
                       and as in Lemma \ref{Gl} we have
 \begin{equation}\label{52N}
           \bold M(G)=\frac{1}{n+1}\int_{Z}\big<\nu, p\big>d\sigma.
 \end{equation}

   \end{cor}
                                                  Based on Corollary \ref{Gl2},
                                                  we can derive a counterpart of Theorem \ref{ToO}
                                                  in this setting involving singularities in boundary.
                                                            \begin{thm}\label{currentOmega}\label{ToO2}
          Let $\Omega\subset \mathbb R^{n+1}$ be a bounded $(n+1)$ domain inducing an integral $(n+1)$-current.
        Suppose that 
                        the graph of  $\eta$ over $\p \Om$
                        is an integral $n$-current in $\R^{n+m+2}$
                        where
          $\eta:\p\Omega\rightarrow \mathbb R^{m+1}$ has  image $N$ of Hausdorff dimension $l<n$.
          If 
          $\eta$
          has
          no continuous extension $\bar \eta$ over $\overline\Omega$ (valued in $\mathbb R^{m+1}$)
      with
           image $\bar \eta(\overline\Omega)$  completely contained in the $\epsilon$-neighborhood of $N$
           for some $\epsilon>0$, 
           then
          there exists $R_\eta\in \mathbb R_+$ such that 
          when constant $R\geq R_\eta$ 
          there is no solution of \eqref{ms} for $f=\eta_R$ on $\p \Omega$.      
          \end{thm}

                      \begin{proof}
                     By the assumption on the graph of $\eta$,
                     the graph of $\eta_R$ for all $R\in R$ is an integral $n$-current
                     and the image of $\eta_R$ is of dimension strictly less than $n$.
                     As a result,
                     every solution $F_i$ of \eqref{ms}
                      for boundary map $\eta_{R_i}$
                      must have
                      an integral $(n+1)$-current $G_i$ (of multiplicity one) as image.
                       
                       Therefore, we can apply Corollary \ref{Gl2} for an upper bound of $\bold M(G_i)$
                       and meanwhile Lemma \ref{dmono} for a lower bound.
                       So argument by contradiction as in the proof of Theorem \ref{ToO}
                       at the beginning of this section
                       leads to the conclusion here.
                      \end{proof}

\begin{figure}[h]
\,\,\,\,\,\,\,\,\includegraphics[scale=0.3]{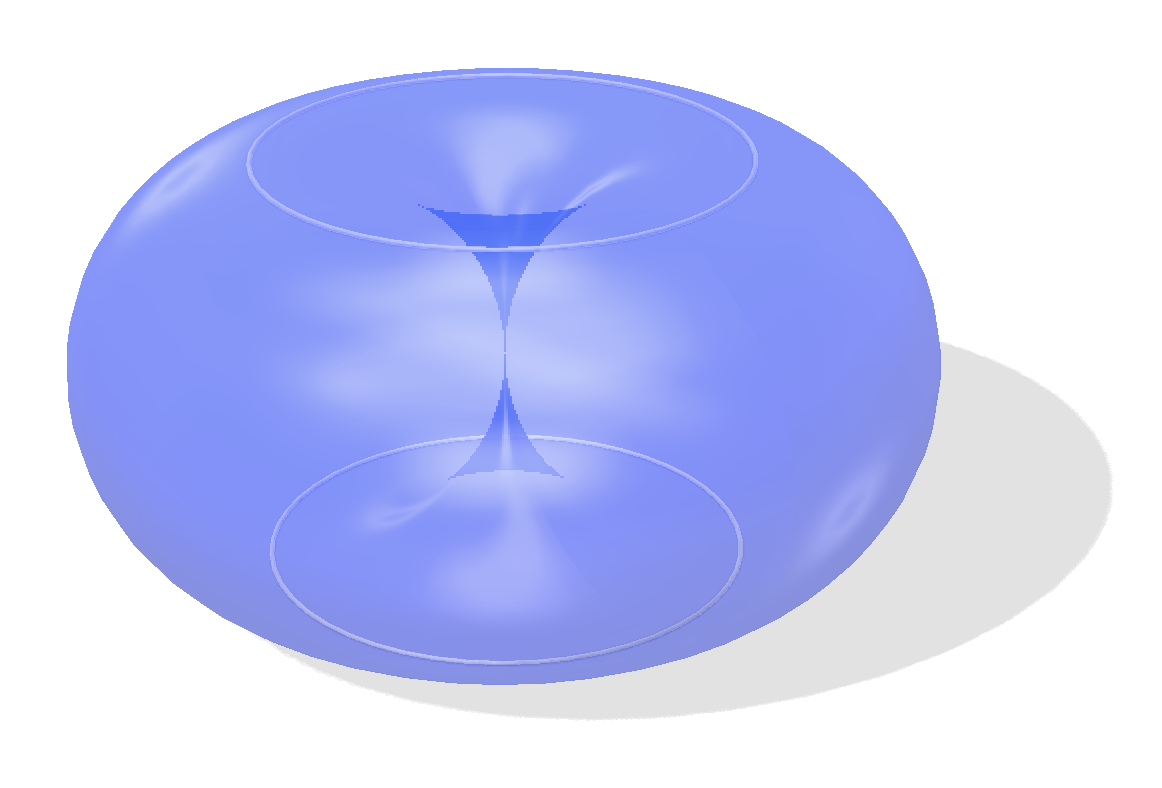}\caption{A picture of $\bold {\check T}^{n+1}$}
\end{figure}

         {\setlength{\parindent}{0cm}
          {\bf Example 4.}
          Take 
          $\Omega=\bold {\check T}^{n+1}=\mathbb D^n(\frac{1}{2})\times S^1(\frac{1}{2})\subset \mathbb R^{n+1}$ 
           with $\eta$ as in Example 2 above
           (now the trivialization of normal bundle along $\{0\}\times S^1(\frac{1}{2})$ 
                          in $\mathbb R\oplus \mathbb R^2$ is given 
                          by $\varepsilon\oplus \mathtt n$
                          where $\mathtt n$ means the trivial normal bundle of $S^1(\frac{1}{2})$ in $\mathbb R^2$).
                          Then Theorem \ref{ToO2} applies for the non-existence of solutions to \eqref{ms}
                          when $f=\eta_R$ for sufficiently large $R$.
          }
                                   \begin{figure}[h]
                              \centering
	\begin{subfigure}[t]{0.4\textwidth}
		\centering
		\includegraphics[scale=0.43]{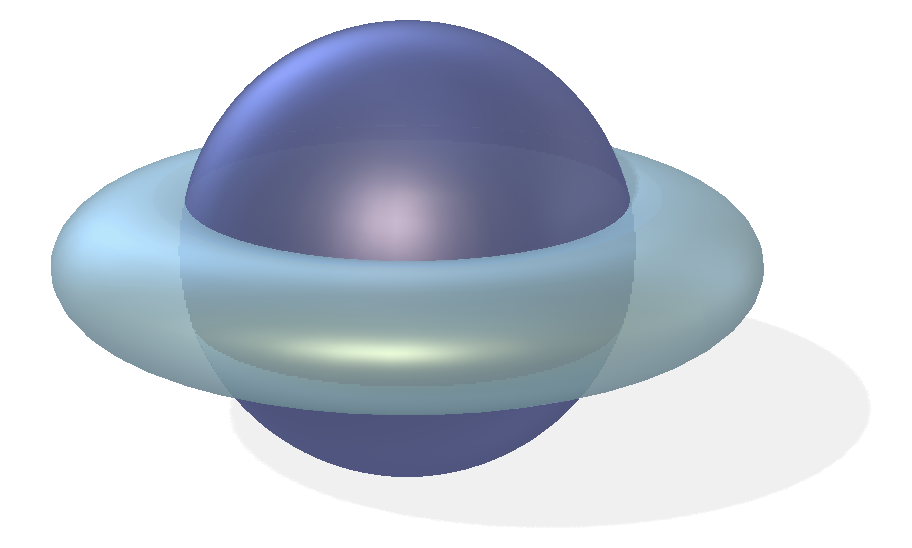}
                              \captionsetup{font={scriptsize}} 
                               \caption{A UFO-style domain}
                               \label{fig:1a}
	\end{subfigure}
	\quad\quad \quad\quad 
	\begin{subfigure}[t]{0.43\textwidth}
		\centering
	 \includegraphics[scale=0.27]{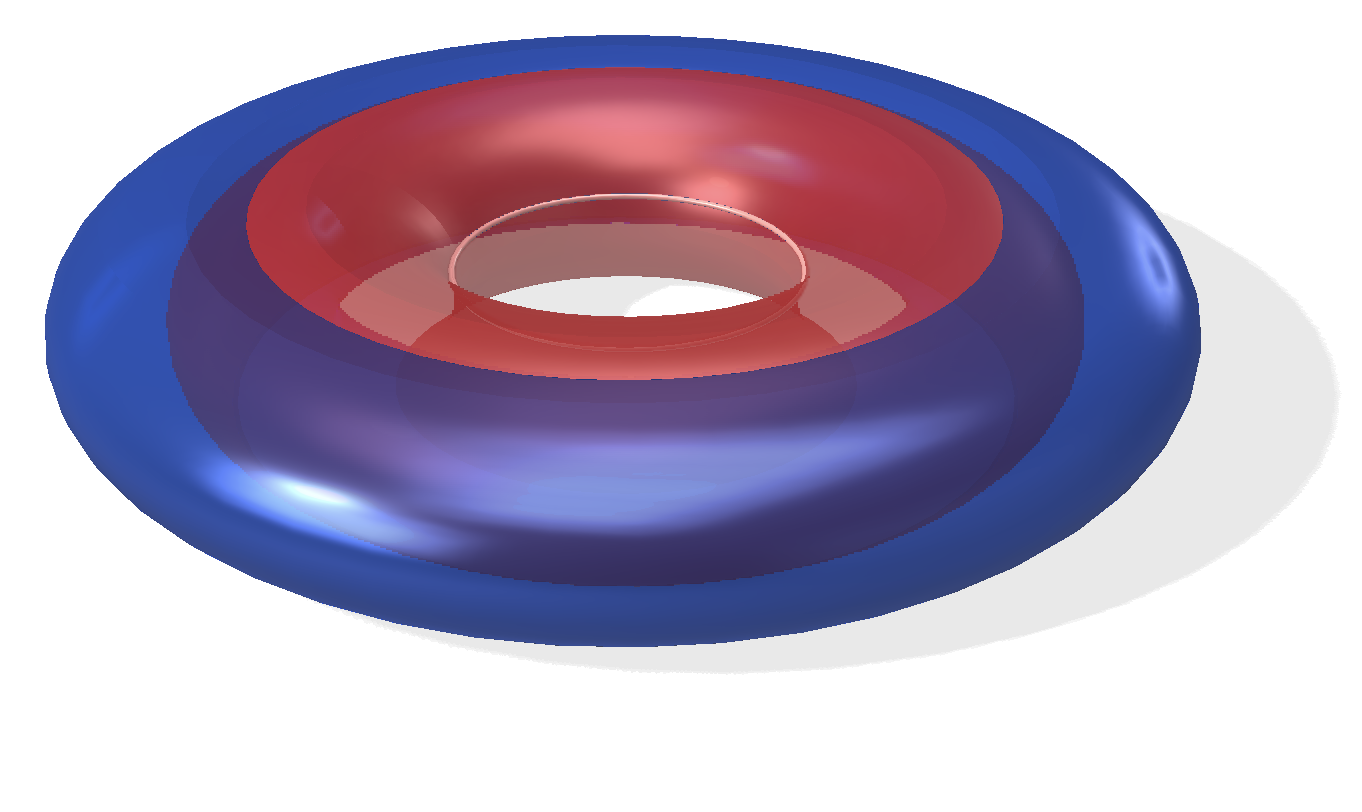}
                           \captionsetup{font={scriptsize}} 
                               \caption{A ``doubled union"}\label{fig:1b}
	\end{subfigure}
                                                           \end{figure}         
                                                           
                   There are also many other non-$C^1$ examples of domains
                   which can be coupled with suitable boundary maps
                   similar to those
                            exhibited in the previous of various homotopy or homology properties
                            for the non-existence type results,
such as
                         a UFO-style domain $\bold {\check T}^{n+1}\bigcup \mathbb D^{n+1}(\frac{1}{2})$, 
                                 a ``doubled union" $\big(1.5\cdot\bold {T}^{n+1}\big)\bigcup\bold {T}^{n+1}$, 
                                 all reasonable polyhedra and et cetera.

              On the other hand, suppose that $\p \Omega$ is of $C^2$ but $\eta$ is $C^2$ except in a Hausdorff-$(n-1)$ measure set $\mathscr S$.
              One can also derive the same kind of non-existence under required boundary conditions.
              However, if a solution  of \eqref{ms}   exists for some $R\cdot\eta$ where $R\in \mathbb R_+$,
              its graph may have tangent cones different from collections of half $(n+1)$-planes (cf. Appendix)
              at points involving $\mathscr S$.
              For instance, take 
              $$\Omega=\mathbb D^{n+1}_{\left(\frac{1}{2},\, 0\cdots\, 0\right)}\left(\dfrac{1}{2}\right)\subset \mathbb D^{n+1}\big(1\big)$$ 
              and 
              boundary data 
              induced by restriction of \eqref{Sg00}
              associated to 
               the (uncountably many) LOMSEs \eqref{Sg1} discovered in \cite{x-y-z0}.
              Now $\mathscr S=\{o\}\subset \p \Omega$,
              and it provides a non-flat tangent cone of graph type at $o$ in $\mathbb R^{m+n+2}$.


{\ }

   \Section{Some connection with complex analysis}{Some connection with complex analysis}\label{S7}

          Note that a complex variety in a K\"ahler manifold is always minimal.
          Let $F:\mathbb D^{2n}(1)\subset \mathbb C^n\rightarrow \mathbb C^m$ be a holomorphic map and $\eta=F|_{\mathbb S^{2n-1}}$ 
          (or other domains).
          Then for any  $R\in\mathbb R$ the Dirichlet problem \eqref{ms} for $f=\eta_R$ has solution $R\cdot F$.  
         
         By Theorem \ref{t2}, we have  the following.
         
                           \begin{cor}\label{c2}
          Let  $F: \mathbb C^n \rightarrow \mathbb C^m$
          be holomorphic and $f=F|_{\mathbb S^{2n-1}}$.
          If the image $N$ of $f$ is an embedded submanifold of dimension $< 2n-1$,
          then $[f]=0\in\pi_{2n-1}(N)$.
          \end{cor}
          
          \begin{rem}
          When $0<m<n$,  the dimension condition on the image of $f$ automatically holds.
          Now, if $[f]\neq 0\in \pi_{2n-1}(N)$, $N$ cannot be an embedded submanifold.
          \end{rem}
          
         By Theorem \ref{Sepa}, we gain another corollary.
         \begin{cor}\label{C12}
         Assume that $n\geq 2$
         and
         $F:\mathbb C^n\rightarrow \mathbb C^m$ is a holomorphic map 
         on $\Om=\mathbb D^{2n}(2)\sim \overline{\mathbb D^{2n}(1)}$
         with $F(S^{2n-1}(1))$ and $F(S^{2n-1}(2))$ of Hausdorff dimensions strictly less than $2n-1$.
         Then it follows $F(S^{2n-1}(1))\bigcap F(S^{2n-1}(2))\neq \emptyset$.
         \end{cor}
         \begin{rem}
         By the Hartogs's phenomenon (e.g. see \cite{g-h}),
         once there is a holomorphic solution $F$
         then its restriction to  $S^{2n-1}(1)$
         is completely decided by that on $S^{2n-1}(2)$.
         So this setting is much  more rigid than that in Remark \ref{notall}.
         \end{rem}

         

{\ }

\appendix
\section*{Appendix}

Lemma \ref{TA} lays the foundation of discussions in this paper.
In this appendix we provide detailed explanations about it
under the assumption that the graph of boundary map $f$ over $\p \Om$ is of $C^{1,\alpha}$.
Obviously, $\p \Om$ needs to be of $C^{1,\alpha}$.

\begin{pfa}
  \begin{figure}[h]
  \,\,\,\,\,\,\,\,\,\,\includegraphics[scale=0.7]{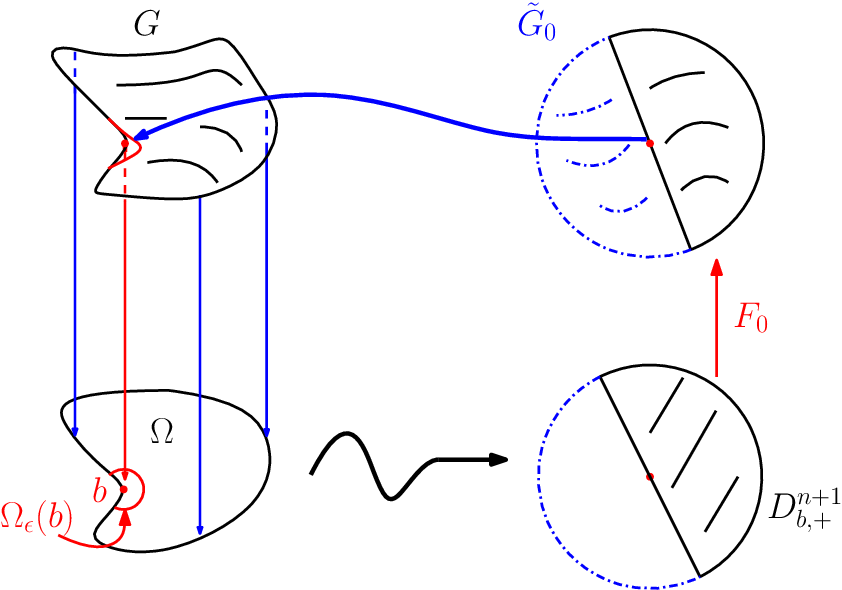}
  \caption{Zoom-in, reflection and density}\label{appfg}
  \end{figure}
      For a point $b\in \p\Omega$,
      consider $\Omega_b(\epsilon):=\mathbb D^{n+1}_b(\epsilon)\bigcap \overline{\Omega}$.
       Let $G$ denote the minimal graph of $F$.
      To blow up $G$ at $(b, F(b))$, we use the size of first $(n+1)$-components instead of Euclidean distance in $\mathbb R^{m+n+2}$.
      Namely, rescale the shifted graph $G-(b, F(b))$ over $\Omega_b(\epsilon)$ by $\frac{1}{\epsilon}$   
      to obtain $G_\epsilon$ the graph of 
           corresponding function $F_\epsilon$ over $\Omega_\epsilon:=\frac{1}{\epsilon}\Omega_b(\epsilon)$ $\subset$ 
      $\mathbb D^{n+1}(1)$.
      As $\epsilon \downarrow 0$,
       $$\Omega_\epsilon\longrightarrow \mathbb D^{n+1}_{b,+}:=\mathbb D^{n+1}(1)\bigcap \{x\in \mathbb D^{n+1}| \big<x,\overrightarrow {n_b}\big>\geq 0\}$$
      where $\overrightarrow {n_b}$ is the interior unit normal of $\p \Omega$ at $b$.

      Since $F$ is a Lipschitz map in the interior of $\Omega$ and continuous up to $\p \Omega$, so is the dilation $F_\epsilon$ with the same Lipschitz constant 
      and all $G_\epsilon$ are contained in a compact set in the cylinder $\mathbb D^{n+1}\times \mathbb R^{m+1}$.
      As $F_{\epsilon_i}$ are bounded Lipschitz maps for $\{\epsilon_i\}\downarrow 0$, there exist a subsequence $\{\epsilon_{i'}\}$ and some Lipschitz $F_0$ defined on $\mathbb D^{n+1}_{b, +}$
      such that
      $F_{\epsilon_{i'}}\rightrightarrows F_0$ on every $K\Subset \mathbb D^{n+1}_{b,+}$.
      \footnote
      {One can fix a genuine  increasing sequence $K_1\subset K_2\subset \cdots \subset K_l\subset \cdots $ of compact sets in $\mathbb D^{n+1}_{b,+}$
      with $\bigcup \mathring{K_l}=\mathbb D^{n+1}_{b,+}$.
      For $K_1$ choose a subsequence of $\{F_{\epsilon_i}\}$ with convergent limit in $K_1$.
      Then for $K_2$ select a further subsequence with convergent limit in $K_2$.
      Repeat the selective procedure by induction for $l\uparrow \infty$ and
       pick up the diagonal sequence $\{F_{\epsilon_{i'}}\}$ for the property of convergence.
      }
      This can be done without requiring convexity on $\Omega$.

      Denote the graph of $F_0$ by $G_0$.
      In our case, $G_0$ coincides with limits (of further subsequents of indices $\{\epsilon_{i'}\}$) in the sense of current or varifold by the compactness result of Federer-Fleming 
      or Allard (cf. \S 6.4 of \cite{A}) for $\{G_{\epsilon_{i'}}\}$. 
      Hence the interior of $G_0$ is a stationary graph over the interior part of $\mathbb D^{n+1}_{b, +}$,
      by taking advantage of Remark 4.11 of \cite{A} in the language of varifolds:
      \begin{quote}
      Suppose $W$ is an open subset of $\mathbb R^{m+n+2}$ and $\lim_{j\rightarrow \infty}V_j=V$ in $\bold V_{n+1}(W)$.
      Evidently, $\|\delta V\|(U)\leq \liminf_{j\rightarrow \infty}\|\delta V_j\|(U)$
      wherever $U$ is an open subset of $W$.
      \end{quote}
      Here $\delta V_j$ stands for the generalized mean curvature of $V_j$.
      Note that $\big<p,\overrightarrow {n_b}\big>\geq 0$ for  $p\in G_0$
      and that the inequality is strict for $p$ in the interior of $G_0$.
      So one can apply the reflection principle in \S3.2 (for currents or varifolds of dimension $\geq 2$) of \cite{a} 
      to extend $F_0$ to $\tilde {F_0}$ defined on the entire $\mathbb D^{n+1}(1)$
      with stationary graph $\tilde {G_0}$, see Figure \ref{appfg}.
Note that 
$$\|\delta(V+\theta_\#V)\|\leq \big(\|\delta V\|+\theta_\#\|\delta V\|\big)\Lcorner \Big(\mathbb R^{m+n+2}\sim \big(T_{(b, F(b))}\p G-(b,F(b))\big)\Big)$$
where $\theta$ is the reflection map about the $n$-dimensional subspace $T_{(b, F(b))}\p G-(b,F(b))$.

Since we want to apply Allard boundary regularity theorem, the key point is to know the density at points in $\p G$.
First, by the Lipschitz constant of the fixed solution $F$, the volume of part of $G$ can be controlled by that of its projection image on $\Omega$.
Since $\Omega_\epsilon\rightarrow \mathbb D^{n+1}_{b,+}$, 
it is essential without any loss to focus on the right column of Figure \ref{appfg} for the density, 
i.e., (for non-convex situation) the ignored parts contribute nothing to the density in the limit as $\epsilon\downarrow 0$. 

     A priori it is unclear to the author whether $\tilde {G_0}$ is a cone. Let VarTan$\big(\tilde {G_0},o\big)$ 
     be the collection of all limit $C$ by a sequence of blowing-ups (of the infinitesimals) of $\tilde {G_0}$ at the origin.
          By the virtue of  Lemma \ref{dmono}, 
     $$\Theta\big(\tilde {G_0}, o\big):=\lim_{\epsilon\rightarrow 0} \Theta\big(\tilde {G_0}, o, \epsilon\big)$$ exists
     with value in $[1,\infty)$.
     Therefore, the stationary $\tilde {G_0}$ satisfies the Corollary 5.7 in \S 8 of \cite{LS} at $o$ and each $C\in$ VarTan$\big(\tilde {G_0}, o\big)$ is a cone.
     
          \begin{figure}[h]
  \,\,\,\,\,\,\,\,\,\,\includegraphics[scale=0.35]{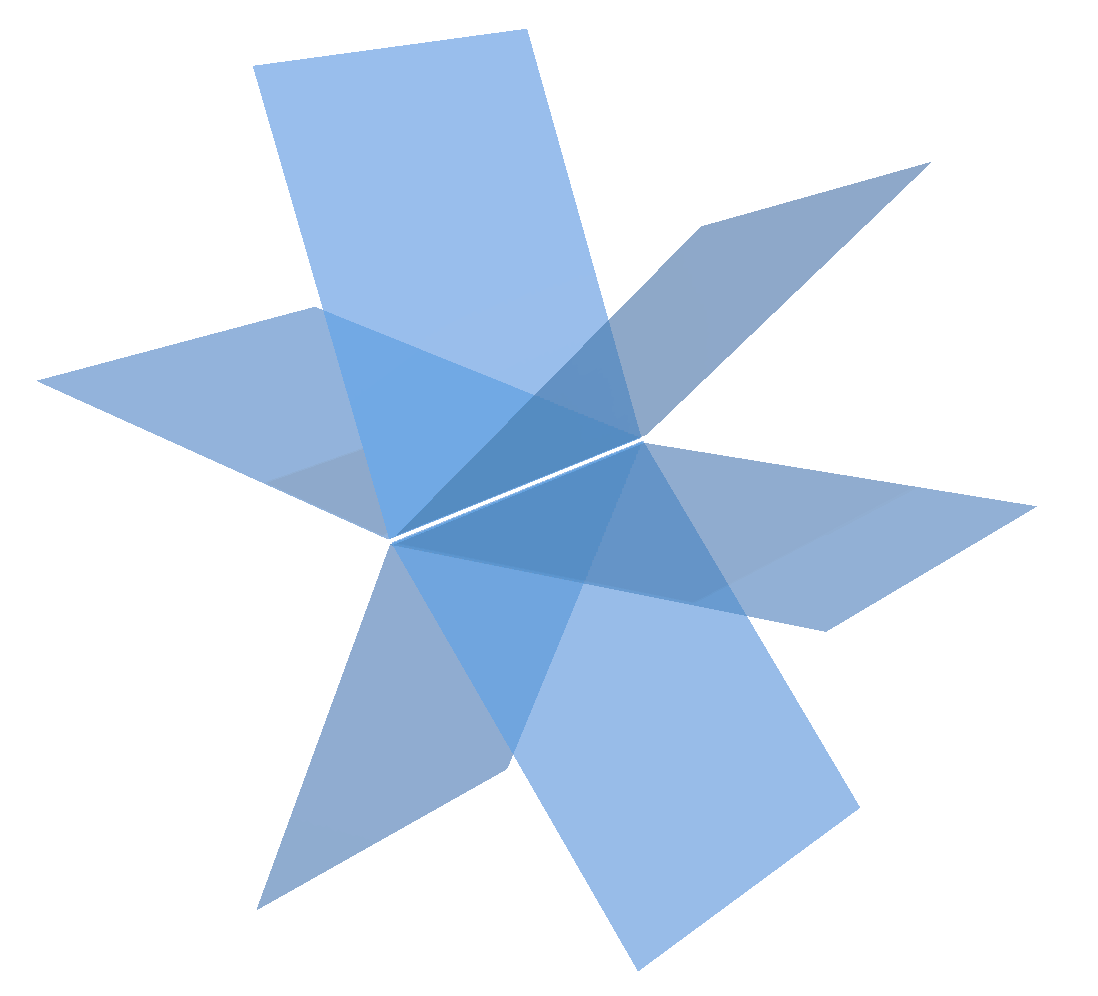}
  \caption{Open Book Structure}\label{OB}
  \end{figure}
     
     Now
     we can apply the beautiful structural result $-$ ``Open Book"  Lemma 5.1 of \cite{a} for such tangent cone to have 
     $$C=\sum_{i=1}^s H_i\,\,\, \,\,\, \text{ as illustrated in Figure }\ref{OB}$$
     where $H_i$ are oriented half $(n+1)$-planes merging in the common edge $T_{(b, F(b))}\p G-(b,F(b))$.
     By property of graph and  reflection principle, $C=H_1+H_2$ and the union of supports of $H_1$ and $H_2$ assembles an $(n+1)$-plane.
     Furthermore, by Theorem 3.5 (1) in \cite{a}, one can apply 3.4 (2) in \cite{A} which asserts that the density of $C$ at $o$ equals double that of $G$ at $(b, F(b))$ (although $G_0$ may not be a cone at this stage) to deduce that
 the density of $G$ at $(b, F(b))$ equals one half.
 
     If we have the $C^2$ assumptions on $\p \Omega$ and boundary maps,
     then
       according to \cite{A, a, Morrey1, Morrey2}
     and the application of Theorem 9.2. of \cite{Morrey0}
the nearby $C^2$ boundary regularity can be set up.

     If we have the $C^{1,\alpha}$ assumptions on $\p \Omega$ and boundary maps,
     then by \cite{B, a2}
     the nearby $C^{1, \alpha}$ boundary regularity can be obtained.
\end{pfa}

\begin{rem}
Note that it merely needs boundary data to be $C^1$ to get the $\frac{1}{2}$ density in the proof.
However, it is unknown to the author whether $C^1$ is enough for a nearby $C^1$ boundary regularity in this minimal graph setting \eqref{ms}.
\end{rem}

{\ }

\bibliographystyle{amsplain}

\end{document}